\documentclass[letterpaper, 10 pt, conference]{ieeeconf}  

\IEEEoverridecommandlockouts                              

\usepackage{amssymb}
\usepackage{amsmath}
\usepackage{bm}

\newtheorem{theorem}{\sc Theorem}[section]
\newtheorem{lemma}[theorem]{\sc Lemma}

\newcommand{\eps}{\varepsilon}

\newcommand{\proofend}{{\medskip\medskip}}

\title{\LARGE \bf
Gaussian Learning-Without-Recall in a Dynamic Social Network
}

\author{Chu Wang$^{1}$ and Bernard Chazelle$^{2}$
\thanks{$^{1}$Nokia Bell Labs, 600-700 Mountain Avenue, Murray Hill, New Jersey 07974, 
{\tt chu.wang@nokia.com }}%
\thanks{$^{2}$Department of Computer Science, Princeton University, 35 Olden Street, Princeton, New Jersey 08540, {\tt chazelle}@{\tt cs.princeton.edu }}%
}

\begin{document}

\maketitle
\thispagestyle{empty}
\pagestyle{empty}

\begin{abstract}

We analyze the dynamics of the Learning-Without-Recall model with Gaussian priors
in a dynamic social network. Agents seeking to learn the state of the world, the ``truth", exchange
signals about their current beliefs across a changing network and update them accordingly. 
The agents are assumed memoryless and rational, meaning that they Bayes-update
their beliefs based on current states and signals, with no other information from the past.
The other assumption is that each agent hears a noisy signal from the truth
at a frequency bounded away from zero.  Under these conditions,
we show that the system reaches truthful consensus almost surely with a convergence rate
that is polynomial in expectation.
Somewhat paradoxically, high outdegree can slow down the learning process.
The lower-bound assumption on the truth-hearing frequency is necessary: even 
infinitely frequent access to the truth offers no guarantee of truthful consensus in the limit.
\end{abstract}

\bigskip

\section{Introduction}

People typically form opinions by updating their current beliefs
and reasons in response to new 
signals from other sources (friends, colleagues, social media, 
newspapers, etc.)~\cite{tahbaz2009learning,acemoglu2011opinion,golub2010naive}.
Suppose there were an information source that made a noisy version
of the ``truth" available to agents connected through a social network.
Under which conditions would the agents reach consensus about their beliefs?
What would ensure truthful consensus (meaning 
that the consensus coincided with the truth)?  How fast would it
take for the process to converge?  To address these questions requires agreeing
on a formal model of distributed learning. 
Fully rational agents update their beliefs by assuming a prior and using Bayes' rule
to integrate all past information available to 
them~\cite{acemoglu2011bayesian,mueller2013general,lobel2015preferences,
mossel2011agreement,banerjee1992simple,bala1998learning}.
Full rationality is intractable in practice~\cite{molavi2015foundations,rahimian2016learning},
so much effort has been devoted to 
developing computationally effective mechanisms, including non- (or partially) Bayesian 
methods~\cite{jadbabaie2012non,molavi2015foundations,
golub2010naive,golub2012homophily,jadbabaie2013information}.
Much of this line of work can be traced back to 
the seminal work of DeGroot~\cite{degroot1974reaching} on linear opinion pooling.

This paper is no exception. Specifically, it follows the {\it Bayesian-Without-Recall} (BWR) model recently 
proposed by Rahimian and Jadbabaie in~\cite{rahimian2016learning};
see also~\cite{rahimian2015learning, rahimian2015log, rahimian2016naive}.
The agents are assumed to be memoryless and rational: this means that they use Bayesian updates 
based on current beliefs and signals with no other information from the past.
The process is local in that agents can collect information only from their neighbors
in a directed graph. In this work, the graph is allowed to change at each time step.
The BWR model seeks to capture the benefits of rational behavior 
while keeping both the computation and the information stored to a minimum~\cite{rahimian2016naive}.

A distinctive feature of our work is that the social network need {\em not} be fixed once and for all.
The ability to modify the communication channels over time reflects the inherently 
changing nature of social networks as well as the reality that our contacts do not all speak to us at once.
Thus even if the underlying network is fixed over long timescales, the model allows
for agents to be influenced by selected subsets of their neighbors.
Dynamic networks are common occurrences in opinion 
dynamics~\cite{hegselmann2002opinion,mohajer2012convergence,chazelle2015, chazelle2015diffusive}
but, to our knowledge, somewhat new in the context of social learning.

Our working model in this paper posits a Gaussian setting:
the likelihoods and initial priors of the agents are normal distributions.
During the learning process, signals are generated as noisy measurements of agents' beliefs
and the noise is assumed normal and unbiased.
Thus all beliefs remain Gaussian at all times~\cite{rahimian2016learning,box2011bayesian}.

Our main result is that,
under the assumption that each agent hears a noisy signal from the truth
at a frequency bounded away from zero, the system reaches truthful consensus
almost surely with a convergence rate
polynomial in expectation.
Specifically, we show that, 
as long as each agent receives a signal from the truth at least once
every $1/\gamma$ steps, the convergence rate is $O(t^{-\gamma/2d})$,
where $d$ is the maximum node outdegree.

Somewhat paradoxically, high outdegree can slow down learning.
The reason is that signals from peer agents are imperfect
conveyors of the truth and can, on occasion, contaminate the network
with erroneous information;
this finding is in line with a similar phenomenon uncovered by
Harel et al.~\cite{harel2014speed}, in which 
social learning system with two Bayesian agents is 
found to be hindered by increased interaction between the agents.
We note that our lower-bound assumption on the truth-hearing frequency is necessary: even 
infinitely frequent access
to the truth is not enough to achieve truthful consensus in the limit.

\medskip
\noindent
\textbf{Further background.} \
Researchers have conducted empirical evaluations of both Bayesian and non-Bayesian 
models~\cite{grimm2014experiment,chandrasekhar2015testing,das2014modeling,bakshy2012role}.
In~\cite{mossel2010efficient}, Mossel et al. analyzed a Bayesian learning system in which each agent 
gets signals from the truth only once at the beginning and then interact with other agents using Gaussian estimators.
In~\cite{moscarini1998social}, Moscarini et al. considered social learning in
a model where the truth is not fixed but is, instead, supplied by a Markov chain.
In the different but related realm of iterated learning, 
agents learn from ancestors and teach descendants.
The goal is to pass on the truth through generations 
while seeking to prevent information loss~\cite{griffiths2007language,smith2009iterated}.

\medskip
\noindent
\textbf{Organization.} \
Section \ref{sec:model} introduces the model and the basic
formulas for single-step belief updates.
Section \ref{sec:mean} investigates the dynamics of the beliefs in expectation
and derive the polynomial upper bound on the convergence rate under the assumption
that each agent hears a signal from the truth at a frequency bounded away from zero.
We demonstrate the necessity of this assumption in~Section \ref{sec:real} 
and prove that the convergence occurs almost surely.

\bigskip

\section{Preliminaries}\label{sec:model}

\subsection{The Model}

We choose the real line $\mathbb{R}$ as the state space
and we denote the agents by $1,2,\ldots, n$;
for convenience, we add an extra agent, labeled $0$, whose belief is a fixed number,
unknown to others, called the {\em truth}.
At time $t=0,1,\ldots$, the belief of agent $i$ 
is a probability distribution over the state space $\mathbb{R}$, which is denoted by $\mu_{t,i}$.
We assume that the initial belief $\mu_{0,i}$ of agent $i$ is Gaussian:
$\mu_{0,i}\sim \mathcal{N}(x_{0,i},\sigma_{0,i}^2)$.
Without loss of generality, we assume the truth is a constant (single-point distribution: $\mu_{t,0}=0$;
$\sigma_{t,0}=0$ for all $t$) and the standard deviation is the same
for all other agents, ie, $\sigma_{0,i}=\sigma_0>0$ for $i>0$.

The interactions between agents are modeled by an infinite sequence $(G_t)_{t\geq 0}$,
where each $G_t$ is a directed graph over the node set $\{0,\ldots, n\}$.
An edge pointing from $i$ to $j$ in $G_t$ indicates that $i$ receives data from $j$ at time $t$.
Typically, the sequence of graphs is specified ahead of time or it is chosen randomly:
the only condition that matters is that it should be independent of the
randomness used in the learning process; specifically, taking expectations and variances
of the random variables that govern the dynamics will assume a fixed graph sequence (possibly random).
Because agent $0$ holds the truth, no edge points away from it.
The adjacency matrix of $G_t$ is denoted by $A_t$: it is an $(n+1)\times(n+1)$ matrix 
whose first row is $(1,0,\dots,0)$.

\subsection{Information Transfer}

At time $t\geq 0$, each agent $i>0$ samples a state $\theta_{t,i}\in {\mathbb R}$ 
consistent with her own belief:  $\theta_{t,i}\sim\mu_{t,i}$.
A noisy measurement $a_{t,i}=\theta_{t,i}+\eps_{t,i}$ is then sent to each
agent $j$ such that $(A_t)_{ji}=1$.
All the noise terms $\eps_{t,i}$ are sampled {\em iid} from $\mathcal{N}(0,\sigma^2)$.
An equivalent formulation is to say that the
likelihood function $l(a|\theta)$ is drawn from $\mathcal{N}(\theta,\sigma^2)$.
In our setting, agent $i$ sends the same data to all of her neighbors; this is done
for notational convenience and the same results would still hold if we were to
resample independently for each neighbor.
Except for the omission of explicit utilities and actions, 
our setting is easily identified as a variant of the BWR model~\cite{rahimian2016learning}.

\subsection{Updating Beliefs}

A single-step update for agent $i>0$ consists of setting $\mu_{t+1,i}$
as the posterior 
${\mathbb P}[\mu_{t,i} | d]\propto {\mathbb P}[d| \mu_{t,i} ]{\mathbb P}[\mu_{t,i}]$,
where $d$ is the data from the neighbors of $i$ received at time $t$.
Plugging in the corresponding Gaussians gives us the classical 
update rules from Bayesian inference~\cite{box2011bayesian}.
Updated beliefs remain Gaussian so we can use the notation
$\mu_{t,i}\sim {\mathcal N}(x_{t,i}, \tau_{t,i}^{-1})$, where $\tau_{t,i}$
denotes the precision $\sigma_{t,i}^{-2}$.
Writing $\tau= \sigma^{-2}$ and letting $d_{t,i}$
denote the outdegree of $i$ in $G_t$, for any $i>0$ and $t\geq 0$, 
\begin{equation}\label{eq:single}
\begin{cases}
\begin{split}
x_{t+1,i} &=(\tau_{t,i}  x_{t,i}+\tau a_1+\dots+\tau a_{d_{t,i}})/(\tau_{t,i}+ d_{t,i} \tau) ; \\
\tau_{t+1,i}&= \tau_{t,i} + d_{t,i}\tau,
\end{split}
\end{cases}
\end{equation}
where $a_1, \ldots, a_{d_{t,i}}$ are the signals received by agent $i$
from its neighbors at time $t$.

\subsection{Expressing the Dynamics in Matrix Form}

Let $D_t$ and $P_t$ denote the $(n+1)$-by-$(n+1)$ diagonal matrices $\text{diag} (d_{t,i})$
and $(\tau_0/\tau)I+ \sum_{k=0}^{t-1}D_k$, respectively, where $I$ is the identity matrix
and the sum is $0$ for $t=0$.
It follows from \eqref{eq:single} that $\mu_{t,i}\sim \mathcal{N}(x_{t,i}, (\tau P_t)_{ii}^{-1})$ for $i>0$.
Regrouping the means in vector form, $\bm{x}_t \! :=(x_{t,0},\dots,x_{t,n})^T$, where
$x_{t,0}=0$ and $x_{0,1},\ldots, x_{0,n}$ are given as inputs, we have
\begin{equation}\label{eq:noisy}
\bm{x}_{t+1}=\left( P_t+D_t\right)^{-1}\left( P_t\bm{x}_t+A_t\left(\bm{x}_t+\bm{u}_t+\bm{\eps}_t\right)\right),
\end{equation}
where $\bm{u}_t$ is such that
$u_{t,0}\sim {\mathcal N}(0,0)$ and, for $i>0$,
$u_{t,i}\sim {\mathcal N}(\bm{0},  (\tau  (P_t)_{ii})^{-1})$;
and $\bm{\eps}_t$ is such that 
$\eps_{t,0}\sim {\mathcal N}(0,0)$ and, for $i>0$,
$\eps_{t,i}\sim {\mathcal N}(\bm{0},  1/\tau)$.
We refer to the vectors $\bm{x}_t$ and $\bm{y}_t \! := \mathbb{E}\, \bm{x}_t$
as the {\em mean process} and the {\em expected mean process}, respectively.
Taking expectations on both sides of~\eqref{eq:noisy}
with respect to the random vectors $\bm{u}_t$ and $\bm{\eps}_t$
yields the update rule for the expected mean process:
$\bm{y}_{0}= \bm{x}_{0}$ and, for $t>0$,
\begin{equation}\label{eq:mean}
\bm{y}_{t+1}=\left( P_t+D_t\right)^{-1}\left( P_t+A_t\right)\bm{y}_t.
\end{equation}
A key observation is that $\left( P_t+D_t\right)^{-1}\left( P_t+A_t\right)$ is a stochastic matrix,
so the expected mean process $\bm{y}_t$ forms a diffusive influence system~\cite{chazelle2015diffusive}:
the vector evolves by taking convex combinations of its own coordinates.  What makes
the analysis different from standard multiagent agreement systems is
that the weights vary over time. In fact, some weights typically tend to 0, which violates
one of the cardinal assumptions used in the analysis of averaging 
systems~\cite{chazelle2015diffusive,moreau05}. This leads us to the use of 
arguments, such as fourth-order moment bounds, that are not commonly
encountered in this area.

\subsection{Our Results}

The belief vector $\bm{\mu}_t$ is Gaussian with mean $\bm{x}_t$ 
and covariance matrix $\Sigma_t$ formed by zeroing out
the top-left element of $(\tau P_t)^{-1}$. 
We say that the system reaches {\em truthful consensus}
if both the mean process $\bm{x}_t$ 
and the covariance matrix tend to zero as $t$ goes to infinity.
This indicates that all the agents' beliefs share a common mean equal to the truth 
and the ``error bars" vanish over time. In view of~\eqref{eq:single},
the covariance matrix indeed tends to $0$ as long as the degrees are nonzero infinitely often,
a trivial condition.  To establish truthful consensus, therefore, boils down to studying
the mean process $\bm{x}_t$. We do this in two parts: first, we show that the expected mean process
converges to the truth; then we prove that fluctuations around it eventually vanish
almost surely.\footnote{The Kullback-Leibler divergence~\cite{jadbabaie2012non} is not suitable here 
because the estimator is Gaussian, hence continuous,
whereas the truth is a single-point distribution.}

\bigskip

\noindent
{\em Truth-hearing assumption}:
Given any interval of length $\kappa \! := \lfloor 1/\gamma \rfloor$, 
every agent $i>0$ has an edge $(i,0)$ in $G_t$ for at least one value of $t$ in that interval.
\bigskip

\begin{theorem}\label{MainTheorem}
Under the truth-hearing assumption,
the system reaches truthful consensus with
a convergence rate bounded by $O(t^{-\gamma/2d})$, where
$d$ is the maximum outdegree over all the networks.
\end{theorem}

\bigskip
We prove the theorem in the next two sections. It will follow
directly from Lemmas~\ref{lemma:y} and~\ref{lemma:x} below.
The convergence rate can be improved to the order of
$t^{-(1-\eps)\gamma /d}$, for arbitrarily small $\eps>0$.
The inverse dependency on $\gamma$ is not surprising: the more access to the truth
the stronger the attraction to it. On the other hand, it might seem
counterintuitive that a larger outdegree should slow down convergence. 
This illustrates the risk of groupthink. It pays to follow the crowds when
the crowds are right. When they are not, however, this distracts from
the lonely voice that happens to be right.

How essential is the truth-hearing assumption? We show that it is necessary.
Simply having access to the truth infinitely often is not enough to achieve truthful consensus.

\bigskip

\subsection{Useful Matrix Inequalities}\label{sec:matrix}

We highlight certain matrix inequalities to be used throughout.
We use the standard element-wise notation
$R\le S$ to indicate that $R_{ij}\le S_{ij}$ for all $i,j$.
The infinity norm $\|R\|_\infty \!  =\max_{i}\sum_{j}|r_{ij}|$
is submultiplicative: $\|RS\|_\infty \! \le\|R\|_\infty\|S\|_\infty$, for any matching
rectangular matrices. On the other hand, the
max-norm $\|R\|_\text{max} \! :=\max_{i,j}|r_{ij}|$ is not, but it is 
transpose-invariant and also satisfies: $\|RS\|_\text{max}\le\|R\|_\infty\|S\|_\text{max}$.
It follows that  
\begin{equation}\label{MatrixIneq}
\begin{split}
\|RSR^T\|_\text{max} \! 
&\le   \|R\|_\infty\|SR^T\|_\text{max}=\|R\|_\infty\|RS^T\|_\text{max}\\
&\le \|R\|_\infty^2\|S^T\|_\text{max}=\|R\|_\infty^2\|S\|_\text{max} .
\end{split}
\end{equation} 

\bigskip

\section{The Expected Mean Process Dynamics}\label{sec:mean}

We analyze the convergence of the mean process in expectation.
The expected mean $\bm{y}_t=\mathbb{E}\, \bm{x}_t$ evolves
through an averaging process entirely determined
by the initial value 
$\bm{y}_0= (0, x_{0,1},\ldots, x_{0,n})^T$ and the graph sequence $G_t$.
Intuitively, if an agent communicates repeatedly with a holder of the
truth, the weight of the latter should accumulate and increasingy influence
the belief of the agent in question.
Our goal in this section is to prove the following result:

\bigskip
\begin{lemma}\label{lemma:y}
Under the truth-hearing assumption,
the expected mean process $\bm{y}_t$ converges to the truth asymptotically.
If, at each step, no agent receives information from more than $d$ agents, 
then the convergence rate is bounded by $C t^{-\gamma/2d}$, where
$C$ is a constant that depends on $\bm{x}_0, \gamma, d, \sigma_0/\sigma$.
\end{lemma}

\bigskip
\noindent{\em Proof.}
We define $B_t$ as the matrix formed by removing the first row and the first column 
from the stochastic $P_{t+1}^{-1}\left(P_t+A_t\right)$.
If we write $\bm{y}_t$ as $(0, \bm{z}_t)$ then, 
by~\eqref{eq:mean},
\begin{equation}\label{eq:mean_matrix}
\begin{pmatrix}
0\\ \bm{z}_{t+1}
\end{pmatrix}=
\begin{pmatrix}
1&\bm{0}\\
\bm{\alpha}_t&B_t
\end{pmatrix}
\begin{pmatrix}
0\\ \bm{z}_{t}
\end{pmatrix},
\end{equation}
where $\bm{\alpha}_{t,i}=(P_{t+1}^{-1})_{ii}$ 
if there is an edge $(i,0)$ at time $t$ and $\bm{\alpha}_{t,i}=0$ otherwise.
This further simplifies to 
\begin{equation}\label{eq:By}
\bm{z}_{t+1}=B_t\bm{z}_t.
\end{equation}
Let $\bm{1}$ be the all-one column vector of length $n$.
Since $P_{t+1}^{-1}\left(P_t+A_t\right)$ is stochastic,
\begin{equation}\label{eq:alpha_B}
\bm{\alpha}_t+B_t\bm{1}=\bm{1}
\end{equation}
In matrix terms, the truth-hearing assumption means that, for any $t\geq 0$,
\begin{equation}\label{eq:alpha}
\bm{\alpha}_{t}+\bm{\alpha}_{t+1}+\dots+\bm{\alpha}_{t+\kappa -1}
\geq Q^{-1}_{t+ \kappa} {\mathbf 1},
\end{equation}
where $Q_t$ is the matrix derived from $P_t$ by removing 
the first row and the last column; the inequality relies on the fact that
$P_t$ is monotonically nondecreasing. 
For any $t>s\ge 0$, we define the product matrix $B_{t:s}$ defined as
\begin{equation}
B_{t:s}:=B_{t-1}B_{t-2}\dots B_{s},
\end{equation}
with $B_{t:t}= I$.
By~\eqref{eq:By},  for any $t>s\geq 0$,
\begin{equation}\label{ztBt}
\bm{z}_t=B_{t:s} \, \bm{z}_s.
\end{equation}
To bound the infinity norm of $B_{t:0}$, we observe that,
for any $0\le l<\kappa-1$, 
the $i$-th diagonal element of $B_{s+\kappa:s+l+1}$
is lower-bounded by
\begin{eqnarray}\label{eq:B_diagonal}
&&  \prod_{j=l+1}^{\kappa-1}(B_{s+j})_{ii}
= \prod_{j=l+1}^{\kappa-1} \frac{(P_{s+j} + A_{s+j})_{ii}}{(P_{s+j+1})_{ii}} \\
&\ge&\prod_{j=l+1}^{\kappa-1} \frac{(P_{s+j})_{ii}}{(P_{s+j+1})_{ii}} 
= \frac{(P_{s+l+1})_{ii}}{ (P_{s+\kappa})_{ii} } \ge \frac{(P_{s})_{ii}}{(P_{s+\kappa})_{ii}}. \nonumber
\end{eqnarray}
The inequalities follow from the nonnegativity of the entries and
the monotonicity of $(P_t)_{ii}$.
Note that~\eqref{eq:B_diagonal} also holds for
$l=\kappa-1$ since $(B_{s+\kappa:s+\kappa})_{ii}=1$.

Since $P_{t+1}^{-1}\left(P_t+A_t\right)$ is stochastic, the row-sum of $B_t$ does not exceed~1;
therefore, by premultiplymultiplying $B_{s+1},B_{s+1},\dots$ on both sides of \eqref{eq:alpha_B},
we obtain:
\begin{equation}\label{B1}
B_{s+\kappa:s}\bm{1} \leq \bm{1}-\sum_{l=0}^{\kappa-1} B_{s+\kappa:s+l+1}\bm{\alpha}_{s+l}.
\end{equation}
Noting that $\|B_t\|_\infty=\|B_t \bm{1}\|_\infty$ for any $t$, as $B_t$ is non-negative,
we combine~\eqref{eq:alpha}, \eqref{eq:B_diagonal}, and \eqref{B1} together to derive:
\begin{equation}\label{eq:B1}
\|B_{s+\kappa:s}\|_\infty\le1-\min_{i>0}\frac{(P_{s})_{ii}}{(P_{s+\kappa})_{ii}^2}.
\end{equation}
%
%
Let $d \! := \max_{t\ge 0} \max _{1\le i\le n} d_{t,i}$
denote the maximum outdegree in all the networks,
and define $\delta=\min\{\tau_0/\tau, 1\}$.
For any $i>0$ and $s\ge \kappa$,
\begin{equation}\label{eq:ps}
\frac{s\delta}{\kappa}\le (P_s)_{ii}\le ds+ \frac{\tau_0}{\tau} ;
\end{equation}
hence,
\begin{equation}\label{eq:ps2}
\max_i (P_{s+\kappa})_{ii} \le d(s+\kappa)+  \frac{\tau_0}{\tau}.
\end{equation}
It follows that
\begin{equation}
\frac{(P_{s+\kappa})_{ii}- (P_s)_{ii}}{(P_{s+\kappa})_{ii}}
=\frac{\sum_{l=0}^{\kappa-1}d_{s+l, i}}{(P_{s+\kappa})_{ii}}
\le\frac{d\kappa^2\delta^{-1}}{s+\kappa}.
\end{equation}
Thus, we have
\begin{equation}\label{eq:poverp}
\begin{split}
\min_{i>0}\frac{(P_{s})_{ii}}{(P_{s+\kappa})_{ii}}
&=   1- \max_{i>0}\frac{(P_{s+\kappa})_{ii}-(P_{s})_{ii}}
    {(P_{s+\kappa}){ii} }\\
&\ge 1- \frac{d \kappa^2 \delta^{-1}}{s+\kappa}.
\end{split}
\end{equation}
We can replace the upper bound of~\eqref{eq:B1} by
$$
1-   \frac{1}{\max_{i>0} (P_{s+\kappa})_{ii}} 
           \min_{i>0}\frac{(P_{s})_{ii}}{(P_{s+\kappa})_{ii}^2} \, ,
$$
which, together with~\eqref{eq:ps2} and \eqref{eq:poverp} gives us
\begin{equation}\label{BskIneq}
\begin{split}
\|B_{s+\kappa:s}\|_\infty
&\le
1- \frac{1}{d(s+\kappa)+\tau_0/\tau} \left(1-\frac{d\kappa^2\delta^{-1}}{s+\kappa}\right) \\
&\le 1-  \frac{1}{2d\kappa(m+2)}.
\end{split}
\end{equation}
The latter inequality holds as long as $s=m\kappa >0$ and
$$m\ge m^*  \! :=    \frac{2d\kappa}{\delta} + \frac{\tau_0}{d\kappa \tau} .
$$
It follows that, for $m_0\ge m^*$,
\begin{equation}\label{eq:B}
\begin{split}
\|B_{(m_0+m)\kappa:m_0 \kappa}\|_\infty 
&\le \prod_{j=2}^{m+1} \left( 1-\frac{1}{2d \kappa (m_0+j)} \right) \\
&\le \exp\left\{ -\frac{1}{2d \kappa}\sum_{j=2}^{m+1}\frac{1}{m_0+j} \right\}.
\end{split}
\end{equation}
The matrices $B_t$ are sub-stochastic so that
$$\|B_t \, \bm{z}\|_\infty\leq \|B_t\|_\infty \|\bm{z}\|_\infty\leq \|\bm{z}\|_\infty.$$
By~\eqref{ztBt}, for any $t\ge (m_0+m)\kappa$,
$$\bm{z}_t = B_{t:  (m_0+m)\kappa}B_{(m_0+m)\kappa:m_0 \kappa}\, \bm{z}_{m_0},$$
so that, by using standard bounds for the harmonic series,
$\ln(k+1)<1+\frac{1}{2}+\dots+\frac{1}{k}\leq 1+ \ln k$, we find that
\begin{equation*}
\begin{split}
\|\bm{z}_t\|_\infty&\leq 
       \|B_{(m_0+m)\kappa:m_0 \kappa}\, \bm{z}_{m_0}\|_\infty \\
&\le \|B_{(m_0+m)\kappa:m_0 \kappa} \|_\infty\|\bm{z}_{0}\|_\infty\nonumber\\ 
&\le  C t^{-1/(2d\kappa)},
\end{split}
\end{equation*}
where $C>0$ depends on $\bm{z}_0, \kappa, d, \tau_0/\tau$.
We note that the convergence rate can be improved to the order of
$t^{-(1-\eps)\gamma /d}$, for arbitrarily small $\eps>0$, by working a little harder with~\eqref{BskIneq}. 
\hfill $\Box$
\proofend

\bigskip

\section{The Mean Process Dynamics}\label{sec:real}

Recall that $\mu_{t,i}\sim {\mathcal N}(x_{t,i}, \tau_{t,i}^{-1})$, where $\tau_{t,i}$
denotes the precision $\sigma_{t,i}^{-2}$.
A key observation about the updating rule in~\eqref{eq:single} is that
the precision $\tau_{t,i}$ is entirely determined by the graph sequence $G_t$
and is independent of the actual dynamics. Adding to this the connectivity property implied by
the truth-hearing assumption, we find immediately that $\tau_{t,i} \rightarrow\infty$ for any agent $i$.
This ensures that the covariance matrix $\Sigma_t$ tends to $0$ as $t$ goes to infinity, which
satisfies the second criterion for truthful consensus.
The first criterion requires that the mean process $\bm{x}_t$ should converge to
the truth~$\bm{0}$. 
Take the vector $\bm{x}_t-\bm{y}_t$ and remove the first coordinate 
$(\bm{x}_t-\bm{y}_t)_0$ to form the vector $\bm{\Delta}_t\in {\mathbb R}^n$.
Under the truth-hearing assumption, we have seen that
$\bm{y}_t \rightarrow \bm{0}$ (Lemma~\ref{lemma:y}),
so it suffices to prove the following:

\medskip

\begin{lemma}\label{lemma:x}
Under the truth-hearing assumption, the deviation
$\bm{\Delta}_t$ vanishes almost surely.
\end{lemma}

\medskip
\noindent{\em Proof.}
We use a fourth-moment argument. The justification for the high order is technical: it is
necessary to make a certain ``deviation power" series converge.
By~\eqref{eq:noisy}, $\bm{x}_t$ is a linear combination of independent Gaussian random vectors $\bm{u}_s$ and $\bm{\eps}_s$ for $0\le s\le t-1$,
and thus $\bm{x}_t$ itself is a Gaussian random vector.
Therefore $\bm{\Delta}_t$ is also Gaussian and its mean is zero.
From Markov's inequality, for any $c>0$,
\begin{equation}\label{eq:markov}
\sum_{t\ge 0}\mathbb{P}[|\Delta_{t,i}|\geq c]\le \sum_{t\ge 0} \frac{\mathbb{E}\, \Delta_{t,i}^4}{c^4}.
\end{equation}
If we are able to show the right hand side of~\eqref{eq:markov} is finite for any $c>0$,
then, by the Borel-Cantelli lemma, with probability one,
the event $|\Delta_{t,i}|\geq c$ occurs only a finite number of times, 
and so $\Delta_{t,i}$ goes to zero almost surely.
Therefore, we only need to analyze the order of the fourth moment $\mathbb{E} \, \Delta_{t,i}^4$.
By subtracting \eqref{eq:mean} from \eqref{eq:noisy}, we have:
\begin{equation}\label{eq:z}
\bm{\Delta}_{t+1}=B_t\bm{\Delta}_t+M_t\bm{v}_t,
\end{equation}
where $\bm{v}_t:=\bm{u}_t+\bm{\eps}_t$ and $M_t:=P_{t+1}^{-1}A_t$;
actually, for dimensions to match, we remove the top coordinate of $\bm{v}_t$ and
the first row and first column of $M_t$ (see previous section for definition of $B_t$).
Transforming the previous identity into a telescoping sum, it follows from
$\bm{\Delta}_0=\bm{x}_0-\bm{y}_0=\bm{0}$ and the definition $B_{t:s}=B_{t-1}B_{t-2}\dots B_s$
that
\begin{equation}\label{eq:z2}
\bm{\Delta}_{t}=\sum_{s=0}^{t-1}B_{t:s+1}M_s\bm{v}_s=\sum_{s=0}^{t-1}R_{t,s}\bm{v}_s,
\end{equation}
where $R_{t,s}:=B_{t:s+1}M_s$.
We denote by $C_1, C_2, \dots$ suitably large 
constants (possibly depending on $\kappa, d, n, \tau, \tau_0$).
By~\eqref{eq:ps}, $\|M_s\|_{\infty}\le C_1/(s+1)$ and,
by~\eqref{eq:B}, for sufficiently large $s$,
$$\|B_{t:s+1}\|_{\infty} \le C_2 (s+1)^{\beta}(t+1)^{-\beta},
$$
where $\beta=1/2d\kappa <1$.
Combining the above inequalities, we obtain the following estimate of $R_{t,s}$ as
\begin{equation}\label{eq:Q}
\|R_{t,s}\|_{\infty}\le C_3(s+1)^{-1+\beta}(t+1)^{-\beta}.
\end{equation}

In the remainder of the proof, the power of a vector is understood element-wise.
We use the fact that
$\bm{v}_s$ and $\bm{v}_{s'}$ are independent if $s\neq s'$
and that the expectation of an odd power of an unbiased Gaussian is always zero.
By Cauchy-Schwarz and Jensen's inequalities,
\begin{eqnarray}\label{eq:Ez4}
&& \hspace{-1cm}\mathbb{E}\, 
         \bm{\Delta}_{t}^4=\left(\sum_{s=0}^{t-1}R_{t,s}\bm{v}_s\right)^4\nonumber\\
&=&\sum_{s=0}^{t-1}\mathbb{E}(R_{t,s} \bm{v}_s)^4+\sum_{0\le s\neq s'<t} 3\, \mathbb{E}(R_{t,s} \bm{v}_s)^2\mathbb{E}(R_{t,s'} \bm{v}_{s'})^2\nonumber\\
&\le&
\sum_{s=0}^{t-1}\mathbb{E}(R_{t,s} \bm{v}_s)^4+ 3\left(\sum_{s=0}^{t-1} \mathbb{E}(R_{t,s} \bm{v}_s)^2\right)^2\nonumber\\
&\le&
\sum_{s=0}^{t-1}\mathbb{E}(R_{t,s} \bm{v}_s)^4+ 3t \, \sum_{s=0}^{t-1} 
        \mathbb{E}^2(R_{t,s} \bm{v}_s)^2 \nonumber\\
&\le&(3t+1)\sum_{s=0}^{t-1}\mathbb{E}(R_{t,s} \bm{v}_s)^4.
\end{eqnarray}
Notice that since the variance of $\bm{v}_t= (v_{t,1},\ldots, v_{t,n})^T$ is nonincreasing, 
there exists a constant $C_4$ such that  $\mathbb{E}\, v_{t,i}^4 \le C_4$.
By Jensen's inequality and the
fact that the variables $v_{t,i}$ are independent for different
values of $i$, we have, for any $i,j,k,l$,
$$  |\mathbb{E}\, v_{t,i} v_{t,j} v_{t,k} v_{t,l}|
\leq  \max_k \mathbb{E}\,  v_{t,k}^4.$$
By direct calculation, it then follows that
\begin{eqnarray}\label{eq:4moment}
&& \hspace{-2.5cm} \max_i \mathbb{E}(R_{t,s} \bm{v}_s)_i^4 
= \max_i \mathbb{E} \Bigl(\sum_{j=1}^n (R_{t,s})_{ij} v_{s,j}\Bigr)^4\nonumber\\
&\le& \max_i \Bigl(\sum_{j=1}^n (R_{t,s})_{i,j}\Bigr)^4 
                   \max_k \mathbb{E}\,  v_{s,k}^4 \nonumber\\
&=& \|R_{t,s}\|^4_\infty \max_k \mathbb{E}\,  v_{s,k}^4 \nonumber\\
&\le& C_5(s+1)^{-4+4\beta}(t+1)^{-4\beta}.
\end{eqnarray}
Summing~\eqref{eq:4moment} over $0\le s\le t-1$, we conclude
from~\eqref{eq:Ez4} that
$\mathbb{E}\, \bm{\Delta}_{t}^4\le C_6 t^{-2}$, and thus
\begin{equation}
\sum_{t\ge 0} \mathbb{E}\, \bm{\Delta}_{t}^4\le C_6 \sum_{t\ge 1}t^{-2}\le C_7.
\end{equation}
By the Borel-Cantelli lemma, it follows that
$\bm{\Delta}_t$ vanishes almost surely.
\hfill $\Box$
\proofend

Theorem~\ref{MainTheorem} follows directly from Lemmas~\ref{lemma:y} and~\ref{lemma:x}.
\hfill $\Box$
\proofend

\bigskip

We now show why the truth-hearing assumption is necessary.
We describe a sequence of graphs $G_t$ that allows every agent 
infinite access to the truth and yet does not lead to truthful consensus.
For this, it suffices to ensure that the expected mean process $\bm{y}_t$ does not converge.
Consider a system with two learning agents with priors $\mu_{0,1}$ and 
$\mu_{0,2}$ from the same distribution $\mathcal{N}(2,1)$. 
We have $x_{0,1}=x_{0,2}= y_{0,1}=y_{0,2}= 2$
and, as usual, the truth is assumed to be 0; the noise variance is $\sigma^2 = 1$.
The graph sequence is defined as follows: set $t_1=0$;
for $k=1,2,\ldots$, agent 1 links to the truth agent at time $t_k$
and to agent 2 at times $t_k+1,\ldots, s_k-1$;
then at time $s_k$, agent 2 links to the truth agent, and
then to agent 1 at times $s_k+1,\ldots, t_{k+1}-1$.
The time points $s_k$ and $t_k$ are defined recursively to ensure that 
\begin{equation}\label{yLBs}
y_{s_k,1}\geq 1+2^{-2k+1} \hspace{.5cm}  \text{and}  \hspace{.5cm}  y_{t_k,2} \geq 1+2^{-2k}.
\end{equation}
In this way, the expected mean processes of the two agents alternate
while possibly sliding down toward 1 but never lower. 
The existence of these time points can be proved by induction.
Since $y_{0,2}=2$, the inequality $y_{t_k,2} \geq 1+2^{-2k}$ holds for $k=1$,
so let's assume it holds up to $k>0$. 
The key to the proof is that, by~\eqref{eq:mean}, 
as agent~1 repeatedly links to agent~2, she is pulled arbitrarily close to it.
Indeed, the transition rule gives us 
$$ 
y_{t+1,1} = \frac{ (P_t)_{11} }{ (P_{t+1})_{11} }\,  y_{t,1}
+ \frac{1}{ (P_{t+1})_{11} }\, y_{t,2},
$$
where $(P_{t+1})_{11}= (P_t)_{11}+1$, which implies that
$y_{t,1}$ can be brought arbitrarily close to $y_{t,2}$ while the latter
does not move: this follows from the fact that any product of the form
$\prod_{t>t_a}^{t_b}\frac{t}{t+1}$ tends to $0$ as $t_b$ 
grows.\footnote{We note that the construction shares 
a family resemblance with one used by Moreau~\cite{moreau05}
to show the non-consensual dynamics of certain multiagent averaging systems. The difference here
is that the weights of the averaging change at each step by increasing
the agent's self-confidence.}

It follows that a suitably increasing sequence of $s_k,t_k$ ensures 
the two conditions~\eqref{yLBs}.
The beliefs of the two agents do not converge to the truth even though they
link to the truth agent infinitely often.

\bigskip\bigskip\bigskip

\footnotesize{
\bibliography{refer}

\begin{thebibliography}{10}

\bibitem{tahbaz2009learning}
Alireza Tahbaz-Salehi, Alvaro Sandroni, and Ali Jadbabaie.
\newblock Learning under social influence.
\newblock In {\em Decision and Control, 2009 held jointly with the 2009 28th
  Chinese Control Conference. CDC/CCC 2009. Proceedings of the 48th IEEE
  Conference on}, pages 1513--1519. IEEE, 2009.

\bibitem{acemoglu2011opinion}
Daron Acemoglu and Asuman Ozdaglar.
\newblock Opinion dynamics and learning in social networks.
\newblock {\em Dynamic Games and Applications}, 1(1):3--49, 2011.

\bibitem{golub2010naive}
Benjamin Golub and Matthew~O. Jackson.
\newblock Naive learning in social networks and the wisdom of crowds.
\newblock {\em American Economic Journal: Microeconomics}, 2(1):112--149, 2010.

\bibitem{acemoglu2011bayesian}
Daron Acemoglu, Munther~A Dahleh, Ilan Lobel, and Asuman Ozdaglar.
\newblock Bayesian learning in social networks.
\newblock {\em The Review of Economic Studies}, 78(4):1201--1236, 2011.

\bibitem{mueller2013general}
Manuel Mueller-Frank.
\newblock A general framework for rational learning in social networks.
\newblock {\em Theoretical Economics}, 8(1):1--40, 2013.

\bibitem{lobel2015preferences}
Ilan Lobel and Evan Sadler.
\newblock Preferences, homophily, and social learning.
\newblock {\em Operations Research}, 64(3):564--584, 2015.

\bibitem{mossel2011agreement}
Elchanan Mossel, Allan Sly, and Omer Tamuz.
\newblock From agreement to asymptotic learning.
\newblock {\em Arxiv preprint arXiv}, 1105, 2011.

\bibitem{banerjee1992simple}
Abhijit~V. Banerjee.
\newblock A simple model of herd behavior.
\newblock {\em The Quarterly Journal of Economics}, pages 797--817, 1992.

\bibitem{bala1998learning}
Venkatesh Bala and Sanjeev Goyal.
\newblock Learning from neighbours.
\newblock {\em The review of economic studies}, 65(3):595--621, 1998.

\bibitem{molavi2015foundations}
Pooya Molavi, Alireza Tahbaz-Salehi, and Ali Jadbabaie.
\newblock Foundations of non-bayesian social learning.
\newblock {\em Columbia Business School Research Paper}, 2015.

\bibitem{rahimian2016learning}
Mohammad~Amin Rahimian and Ali Jadbabaie.
\newblock Learning without recall from actions of neighbors.
\newblock In {\em 2016 American Control Conference (ACC)}, pages 1060--1065.
  IEEE, 2016.

\bibitem{jadbabaie2012non}
Ali Jadbabaie, Pooya Molavi, Alvaro Sandroni, and Alireza Tahbaz-Salehi.
\newblock Non-bayesian social learning.
\newblock {\em Games and Economic Behavior}, 76(1):210--225, 2012.

\bibitem{golub2012homophily}
Benjamin Golub and Matthew~O. Jackson.
\newblock How homophily affects the speed of learning and best response
  dynamics.
\newblock 2012.

\bibitem{jadbabaie2013information}
Ali Jadbabaie, Pooya Molavi, and Alireza Tahbaz-Salehi.
\newblock Information heterogeneity and the speed of learning in social
  networks.
\newblock {\em Columbia Business School Research Paper}, (13-28), 2013.

\bibitem{degroot1974reaching}
Morris~H. DeGroot.
\newblock Reaching a consensus.
\newblock {\em Journal of the American Statistical Association},
  69(345):118--121, 1974.

\bibitem{rahimian2015learning}
Mohammad~Amin Rahimian, Shahin Shahrampour, and Ali Jadbabaie.
\newblock Learning without recall by random walks on directed graphs.
\newblock In {\em 2015 54th IEEE Conference on Decision and Control (CDC)},
  pages 5538--5543. IEEE, 2015.

\bibitem{rahimian2015log}
Mohammad~Amin Rahimian et~al.
\newblock Learning without recall: A case for log-linear learning.
\newblock {\em IFAC-PapersOnLine}, 48(22):46--51, 2015.

\bibitem{rahimian2016naive}
Mohammad~Amin Rahimian and Ali Jadbabaie.
\newblock Naive social learning in ising networks.
\newblock In {\em 2016 American Control Conference (ACC)}, pages 1088--1093.
  IEEE, 2016.

\bibitem{hegselmann2002opinion}
Rainer Hegselmann and Ulrich Krause.
\newblock Opinion dynamics and bounded confidence models, analysis, and
  simulation.
\newblock {\em Journal of Artificial Societies and Social Simulation}, 5(3),
  2002.

\bibitem{mohajer2012convergence}
Soheil Mohajer and Behrouz Touri.
\newblock On convergence rate of scalar hegselmann-krause dynamics.
\newblock In {\em Proceedings of the IEEE American Control Conference (ACC)},
  2013.

\bibitem{chazelle2015}
Bernard Chazelle and Chu Wang.
\newblock Inertial {H}egselmann-{K}rause systems.
\newblock In {\em Proceedings of the IEEE American Control Conference (ACC)},
  pages 1936--1941, 2016.

\bibitem{chazelle2015diffusive}
Bernard Chazelle.
\newblock Diffusive influence systems.
\newblock {\em SIAM Journal on Computing}, 44(5):1403--1442, 2015.

\bibitem{box2011bayesian}
George~E.P. Box and George~C. Tiao.
\newblock {\em Bayesian inference in statistical analysis}, volume~40.
\newblock John Wiley \& Sons, 2011.

\bibitem{harel2014speed}
Matan Harel, Elchanan Mossel, Philipp Strack, and Omer Tamuz.
\newblock The speed of social learning.
\newblock {\em arXiv preprint arXiv:1412.7172}, 2014.

\bibitem{grimm2014experiment}
Veronika Grimm and Friederike Mengel.
\newblock An experiment on belief formation in networks.
\newblock {\em Available at SSRN 2361007}, 2014.

\bibitem{chandrasekhar2015testing}
Arun~G. Chandrasekhar, Horacio Larreguy, and Juan~Pablo Xandri.
\newblock Testing models of social learning on networks: Evidence from a lab
  experiment in the field.
\newblock Technical report, National Bureau of Economic Research, 2015.

\bibitem{das2014modeling}
Abhimanyu Das, Sreenivas Gollapudi, and Kamesh Munagala.
\newblock Modeling opinion dynamics in social networks.
\newblock In {\em Proceedings of the 7th ACM international conference on Web
  search and data mining}, pages 403--412. ACM, 2014.

\bibitem{bakshy2012role}
Eytan Bakshy, Itamar Rosenn, Cameron Marlow, and Lada Adamic.
\newblock The role of social networks in information diffusion.
\newblock In {\em Proceedings of the 21st international conference on World
  Wide Web}, pages 519--528. ACM, 2012.

\bibitem{mossel2010efficient}
Elchanan Mossel and Omer Tamuz.
\newblock Efficient bayesian learning in social networks with gaussian
  estimators.
\newblock {\em arXiv preprint arXiv:1002.0747}, 2010.

\bibitem{moscarini1998social}
Giuseppe Moscarini, Marco Ottaviani, and Lones Smith.
\newblock Social learning in a changing world.
\newblock {\em Economic Theory}, 11(3):657--665, 1998.

\bibitem{griffiths2007language}
Thomas~L. Griffiths and Michael~L. Kalish.
\newblock Language evolution by iterated learning with bayesian agents.
\newblock {\em Cognitive Science}, 31(3):441--480, 2007.

\bibitem{smith2009iterated}
Kenny Smith.
\newblock Iterated learning in populations of bayesian agents.
\newblock In {\em Proceedings of the 31st annual conference of the cognitive
  science society}, pages 697--702. Citeseer, 2009.

\bibitem{moreau05}
Luc Moreau.
\newblock Stability of multiagent systems with time-dependent communication
  links.
\newblock {\em IEEE Transactions on Automatic Control}, 50:169--182, 2005.

\end{thebibliography}
\bibliographystyle{unsrt}}

\end{document}